\documentclass[12pt,reqno]{amsart}
\usepackage{amsthm,amsfonts,amssymb,amscd,amsmath,euscript}
\usepackage[all]{xy}
\newtheorem{theorem}[subsection]{Theorem}

\newtheorem{proposition}[subsection]{Proposition}

\newtheorem{conjecture}[subsection]{Conjecture}
\newtheorem{lemma}[subsection]{Lemma}
\newtheorem{corollary}[subsection]{Corollary}

\theoremstyle{definition}

\newtheorem{proposition-definition}[subsection]{Proposition-Definition}

\theoremstyle{remark}
\newtheorem{remark}[subsection]{Remark}

\newcommand{\dual}{{\scriptscriptstyle \vee}}
\newcommand{\sdual}{{\,\check{}}}

\newcommand{\fatdot}{{\scriptscriptstyle \bullet}}

\newcommand{\Sd}{{S^{(d)}}}
\newcommand{\Sdd}{{S^{[d]}}}
\newcommand{\smallcup}{{\raisebox{0.2 ex}{$\scriptscriptstyle \cup$}}}

\newcommand{\Xdd}{{X^{[d]}}}

\newcommand{\Br}{\operatorname{Br}\nolimits}

\newcommand{\ch}{\operatorname{ch}\nolimits}

\newcommand{\Coh}{\mathfrak C\mathfrak o\mathfrak h}
\newcommand{\coh}{{\operatorname{coh}\nolimits}}

\newcommand{\Div}{\operatorname{Div}}

\newcommand{\et}{{\operatorname{\acute et}\nolimits}}

\newcommand{\End}{\operatorname{End}\nolimits}
\newcommand{\Ext}{\operatorname{Ext}\nolimits}
\newcommand{\GCD}{\operatorname{g.c.d.\:}\nolimits}

\newcommand\Hilb{{\operatorname{Hilb}\nolimits}}

\newcommand{\HOM}{{{\EuScript H}om\:}}
\newcommand{\Hom}{\operatorname{Hom}\nolimits}
\newcommand{\id}{\operatorname{id}\nolimits}

\newcommand{\Mod}{\mathfrak M\mathfrak o\mathfrak d}

\newcommand{\Pic}{\operatorname{Pic}\nolimits}

\newcommand{\rk}{\operatorname{rk}\nolimits}

\newcommand{\tor}{{\operatorname{tor}}}

\newcommand{\Supp}{\operatorname{Supp}}

\newcommand{\Td}{\operatorname{Td}\nolimits}

\newcommand{\CC}{{\mathbb C}}

\newcommand{\ZZ}{{\mathbb Z}}
\newcommand{\QQ}{{\mathbb Q}}
\newcommand{\PP}{{\mathbb P}}

\newcommand{\DDB}{{\mathbf D}}

\newcommand{\LD}{{\mathbf L}}
\newcommand{\R}{{\mathbf R}}
\newcommand{\RD}{{\mathbf R}}

\newcommand{\OOO}{{\EuScript O}}

\newcommand{\III}{{\EuScript I}}

\newcommand{\GGG}{{\EuScript G}}
\newcommand{\EEE}{{\EuScript E}}

\newcommand{\LLL}{{\EuScript L}}
\newcommand{\FFF}{{\EuScript F}}

\newcommand{\NNN}{{\EuScript N}}
\newcommand{\MMM}{{\EuScript M}}
\newcommand{\PPP}{{\EuScript P}}

\newcommand{\EXT}{{\EuScript Ext}}

\newcommand{\vv}{{\boldsymbol{ v}}}

\renewcommand{\bar}[1]{\overline{#1}}

\newcommand{\Lotimes}{\stackrel{\mathbf L}{\otimes}}

\newcommand\alp{\alpha}

\newcommand\eps{\epsilon}

\renewcommand\phi{\varphi}

\newcommand{\isoto}{{\lra\hspace{-1.3 em}
\raisebox{ 0.6 ex}{$\textstyle\sim$}\hspace{0.8 em}}}
\newlength{\rrrr}
\newcommand{\into}{\hookrightarrow}

\newcommand\lra{{\longrightarrow}}

\newcommand\rar{\rightarrow}
\newcommand\lrdash{\:
\xymatrix@1{\ar@{-->}[r]&}\:
}
\newcommand{\lrdashar}[1]{\:
\xymatrix@1{\ar@{-->}[r]^{#1}&}\:
}

\renewcommand\emptyset{\varnothing}

\author{D. Markushevich
}
\address{\scriptsize Dimitri Markushevich:  
Math\'ematiques - b\^{a}t.M2, Universit\'e Lille 1, 
F-59655 Villeneuve d'Ascq Cedex, France}
\email{markushe@math.univ-lille1.fr}


\title{Rational Lagrangian fibrations on punctual Hilbert schemes of K3 surfaces}

\begin{document}

\begin{abstract}
A rational Lagrangian fibration $f$ on an irreducible symplectic variety $V$
is a rational map which is birationally equivalent to a regular surjective
morphism with Lagrangian fibers. 
By analogy with K3 surfaces, it is natural to expect that
a rational Lagrangian fibration exists if and only if $V$ has a divisor
$D$ with Bogomolov--Beauville square 0. This conjecture is proved in the case
when $V$ is the Hilbert scheme of $d$ points on a
generic K3 surface $S$ of genus $g$ under the hypothesis that
its degree $2g-2$ is a square times $2d-2$.
The construction of $f$ uses a twisted Fourier--Mukai transform which induces
a birational isomorphism of $V$ with a certain moduli
space of twisted sheaves on another
K3 surface $M$, obtained from $S$ as its Fourier--Mukai partner.
\end{abstract}

\maketitle


\setcounter{section}{-1}

\section{Introduction}

According to Beauville \cite{Beau-1}, \cite{Beau-2}, the $d$-th symmetric
power $\Sd$ of a K3 surface $S$ has a natural resolution of
singularities, the punctual Hilbert scheme $\Sdd=\Hilb^dS$, which
is a $2d$-dimensional irreducible symplectic variety.
Here a holomorphically symplectic manifold is called {\em an irreducible symplectic variety}
if it is projective and simply connected and has a unique
symplectic structure up to proportionality. In dimension 2,
the irreducible symplectic varieties are exactly K3 surfaces.
Here and throughout the paper, a $K3$ {\em surface} means a 
{\em projective K3 surface}, and all the {\em varieties} are
assumed to be projective.

The 2-dimensional cohomology $H^2(V,\ZZ)$ of an irreducible symplectic variety $V$
has an integral quadratic form $q_V$ with remarkable properties, called the Bogomolov--Beauville form.
This form together with the Hodge structure on $H^2(X)$ gives rise to many striking analogies
between K3 surfaces, where $q_V$ is just the intersection form, and 
higher-dimensional irreducible symplectic varieties. They have very similar deformation
theories, descriptions of period maps, Torelli theorems, structures of the ample (or K\"ahler) cone
(see \cite{Bou}, \cite{Hu-1}, \cite{Hu-2}, \cite{O'G-1}, \cite{O'G-2}, \cite{S-1},
\cite{S-2}  and references therein). There are also differences. For example,
there is only one deformation class of K3 surfaces, but at least two deformation classes
of irreducible symplectic varieties of any even dimension $>2$; one of them is the class
of the variety $\Sdd$. Another difference is the existence of nontrivial birational
isomorphisms (flops) between irreducible symplectic varieties. Huybrechts shows in \cite{Hu-0}
that two birationally equivalent irreducible symplectic varieties not only are deformation
equvalent and have the same period, but also represent nonseparated points of the moduli space.

It is natural to expect that the analogy between the K3 surfaces and
irreducible symplectic varieties has many more manifestations that are
still to be proved. For example, O'Grady conjectures in \cite{O'G-2} that
an irreducible symplectic variety $V$, deformation equivalent to $\Sdd$,
has an antisymplectic birational involution if it has
a divisor class with Bogomolov--Beauville square 2. This is of course true
for K3 surfaces. He also raises the problem of finding explicit geometric
constructions for the irreducible symplectic varieties having a divisor class
with small square. 

A more specific question in this direction is the characterization of
varieties $V$ that have a structure of a fibration, that is a regular
surjective map $f:V\rar B$ to some other variety $B$, different from a point, with connected
fibers of positive dimension. If $V$ is K3, then such a map exists
if and only if $V$ has a divisor $D$ with square 0, and then 
the divisor $D'$ obtained from $\pm D$ by a number of reflections in $(-2)$-curves
defines the structure of an elliptic pencil $\phi_{|D'|}:V\rar \PP^1$.

One cannot generalize this straightforwardly by saying that a higher
dimensional irreducible symplectic variety admits a fibration
if and only if the Bogomolov--Beauville form on $\Pic (V)$ represents
zero. Indeed, it follows from the description of $\Pic (\Sdd)$ in
Theorem \ref{H2-for-Xdd} that if $S$ is a generic K3 surface of degree $2d-2$,
then $\Pic (\Sdd)$ has exactly two primitive effective divisor classes $h\pm e$
of square zero and none of them defines a fibration structure. Nevertheless,
$V=\Sdd$ admits a structure of a rational fibration, that is a rational map $f:V\dasharrow B$
which can be birationally transformed to a genuine fibration $g:W\rar B$ on a symplectic
variety $W$ birational to $V$. This map is introduced by formula (\ref{theta}) and it coincides with
the rational map $\phi_{|h-e|}$ given by the complete linear system $|h-e|$.

Thus the expected generalization is the following: an irreducible symplectic variety $V$
has a structure of a {\em rational} fibration if and only if
it has a divisor $D$ with square 0. 
Matsushita \cite{Mat-1}, \cite{Mat-2} has proved important results on the structure of
regular fibrations $f:V\rar B$ on an irreducible symplectic variety $V$: $\dim B=\frac{1}{2}\dim V$
and the generic fiber of $f$ is an abelian variety which is a Lagrangian subvariety of $V$
with respect to the symplectic structure. If, moreover, $B$ is nonsingular, then $B$ has
the Hodge numbers of a projective space. Remark that no examples are known of fibrations
on irreducible symplectic varieties with base different from a projective space.
So one might complete the conjecture in saying that the base $B$ of any rational fibration $f$
on an irreducible symplectic variety $V$ of dimension $2n$ is the projective space $\PP^n$
and that $f$ is given by the complete linear system $|D|$ of a divisor $D$ with 
$q_V(D)=0$.

In the present article, we prove
this conjecture for the varieties $\Sdd$ constructed from {\em generic}
primitively polarized K3 surfaces $S$ of degree $(2d-2)m^2$, $m\geq 2$
(Corollary \ref{hminusme}). In the next few lines we describe briefly our construction.
We define another $K3$ surface $M$, which is a moduli
space of sheaves on $S$, and a birational map $\mu$ from $\Sdd$ to another
irreducible symplectic variety $V$. The latter is a moduli
space of $\alp$-twisted sheaves on $M$ for some element $\alp$ of the Brauer group
$\Br (M)$, and $\mu$ is induced by the twisted Fourier--Mukai transform
defined by a twisted universal sheaf on $S\times M$. Further, $V$
is a torsor under the relative Jacobian of a linear system $|C|$ on $M$, and hence
has a natural morphism \mbox{$f:V\rar |C|\simeq\PP^d$}
which sends each twisted sheaf to its support. This map is a {\em regular}
Lagrangian fibration, and the {\em rational} one on $\Sdd$ is $\pi=f\circ\mu$.

The same result was obtained independently and almost simultaneously in
\cite{S-2}. Later Yoshioka proved the regularity of this rational Lagrangian
fibration; Yoshioka's proof is included
in the version 3 of loc. cit. The paper \cite{Gu} contains a similar result
for the other series of Beauville's examples, generalized Kummer varieties;
in this case the proof does not necessitate a use of twisted sheaves,
and it seems very likely that the rational Lagrangian fibration constructed
by the author is indeed not regular.

By Proposition \ref{H2-for-Xdd}, $\Sdd$ has a divisor class with Bogomolov--Beauville
square zero for a generic primitively polarized K3
surface $S$ of degree $n$ if and only if $k^2n=(2d-2)m^2$
for some integers $k\geq 1$, $m\geq 1$, $d\geq 2$ with relatively prime $k,m$.
Regular Lagrangian fibrations on $\Sdd$ have been known before in some particular
cases when $k=1$. Hassett and Tschinkel \cite{HTsch-1},  \cite{HTsch-2}
have proved the existence of a regular Lagrangian fibration on $\Sdd$ for a generic K3 surface $S$ of
degree $2m^2$, which corresponds to $d=2$, $m\geq 2$. The authors of \cite{IR}
provided an explicit construction of such a fibration in the case $d=3$, $m=2$.
No examples are known with $k>1$.

In Section 1, we gather generalities on irreducible symplectic varieties and
fibrations on them. In Section 2, we cite necessary notions and results on
Fourier--Mukai transforms and on
twisted sheaves following  \cite{Cal-2}, \cite{HS}, \cite{Y-2}. In Section 3 we
define the moduli K3 surface $M=M(m,H,(d-1)m)$ and study the properties of the sheaves
on the initial surface $S$ represented by points of $M$. In Section 4, the main result
(Theorem \ref{main}, Corollary \ref{hminusme}) is proved. In conclusion, we show that the same construction applied to nongeneric K3 surfaces of degree \mbox{$(2d-2)m^2$} yields nonregular
rational fibration maps $f$ (Proposition~\ref{lattice-polar}).

\bigskip

{\sc Acknowledgements.} I thank A.~Iliev for extensive discussions:
in fact, we started a joint work on Lagrangian fibrations, but finally
divided our problem into two parts and continued the work separately.
I~am grateful to D.~Orlov for introducing me to the twisted sheaves.
I acknowledge discussions with K.~Ranestad and S.~Popescu,
as well as the hospitality of the Mathematishches Forschgungsinstitut
in Oberwolfach, where a part of the work was done. I would also
like to thank the referee for a thorough reading of the manuscript
and for suggestions which helped to correct some imprecisions.

\section{Preliminaries}\label{prelim}

A {\em symplectic variety} is a nonsingular
projective variety $V$ over $\CC$ having a nondegenerate
holomorphic 2-form $\alp\in H^0(\Omega^2_V)$. It is called
{\em irreducible symplectic} if it is, moreover, simply connected
and \mbox{$h^{0,2}(V)=1$.} By the Bogomolov--Beauville decomposition
theorem \cite{Beau-2}, every symplectic variety becomes, after a finite
\'etale covering, a product of a complex torus and a number of
irreducible symplectic varieties.

\begin{theorem}[Beauville, \cite{Beau-2}]
Let $V$ be an irreducible symplectic variety of dimension $2d$. Then there
exists a constant $c_V$ and an integral idivisible quadratic form $q_V$ of
signature $(3,b_2(V)-3)$
on the cohomology $H^2(V,\ZZ)$ such that $\gamma^{2d}=c_Vq_V(\gamma)^d$
for all $\gamma\in H^2(V,\ZZ)$, where $\gamma^{2d}\in H^{4d}(V,\ZZ)$
denotes the 2d-th power of $\gamma$ with respect to the cup product
in $H^*(V,\ZZ)$, and $H^{4d}(V,\ZZ)$ is naturally identified with $\ZZ$.
\end{theorem}

The form $q_V$ was first introduced by Bogomolov
in \cite{Bo}, so we will call it
Bogomolov--Beauville form, and $c_V$ is called Fujiki's constant.

In dimension 2, the irreducible symplectic varieties are just K3 surfaces.
Historically, the first constructions of higher-dimensional
irreducible symplectic varieties
belong to Fujiki \cite{F} (one example of dimension $4$)
and Beauville
\cite{Beau-1}, \cite{Beau-2} (two infinite series
of examples in all even dimensions $\geq 4$).
The Beauville's examples are: 1)~$X^{[d]}=\Hilb^d(X)$, the Hilbert scheme of
\mbox{$0$-dimensional} subschemes of length $d$ in a K3 surface $X$, and  2) $K_n(A)$,
the generalized Kummer variety associated to an abelian surface $A$.
The latter is defined as the fiber of the summation map $A^{[n+1]}\rar A$.
The punctual Hilbert scheme $\Xdd$ has a natural Hilbert--Chow map
$\Xdd\rar X^{(d)}$ sending each \mbox{$0$-dimensional} subscheme to
the associated 0-dimensional cycle of degree $d$, considered as
a point of the $d$-th symmetric power $X^{(d)}$ of $X$.
The Hilbert--Chow map is a resolution of singularities whose
exceptional locus is a single irreducible divisor $E$, the inverse
image of the big diagonal of $X^{(d)}$. By a {\em lattice} we mean
a free $\ZZ$-module of finite rank endowed with a nondegenerate integer
quadratic form.

\begin{proposition}\label{H2-for-Xdd} Let $X$ be a K3 surface.
Then $c_{\Xdd}=\frac{(2d)!}{d!2^d}$ and there is a natural isomorphism of lattices
$H^2(\Xdd,\ZZ)\simeq H^2(X,\ZZ)\stackrel{\perp}{\oplus}\ZZ e$, \ $e^2=-2(d-1)$,
where $e$ is the class of the exceptional divisor $E$ of the Hilbert--Chow
resolution, and $e^2$ stands for the square of $e$ with respect to the Bogomolov--Beauville
form $q_{\Xdd}$.
\end{proposition}

\begin{proof}
See \cite{HL}, 6.2.14.
\end{proof}

Using the isomorphism of Proposition \ref{H2-for-Xdd}, we will 
denote a class in $\Pic X$ or $H^2(X,\ZZ)$
and its image in $\Pic \Xdd$ or $H^2(\Xdd,\ZZ)$ by the same symbol.

Other examples of irreducible symplectic varieties are given
by moduli spaces of sheaves on a K3 or abelian surface $Y$.
Mukai \cite{Mu-1} has endowed the integer cohomology
$H^*(Y)$  with the following bilinear form:
\begin{equation}\label{mukai-form}
\langle(v_0,v_1,v_2),(w_0,w_1,w_2)\rangle=v_1\smallcup w_1-
v_0\smallcup w_2-v_2\smallcup w_0,
\end{equation}
where $v_i,w_i\in H^{2i}(Y)$. We will denote $\langle v,v\rangle$
simply by $v^2$. For a sheaf $\FFF$ on $Y$,
the {\em Mukai vector} is 
\begin{multline*}
v(\FFF)=\ch(\FFF)\sqrt{\Td(Y)}=(\rk \FFF, c_1(\FFF),
\chi (\FFF)-\eps\rk \FFF)\in \\
H^0(Y)\oplus H^2(Y)\oplus H^4(Y)=H^*(Y),
\end{multline*} 
where $\Td(Y)$ is the
Todd genus, $H^4(Y)$ is naturally identified
with $\ZZ$ and $\eps=0$, resp. $1$ if $Y$ is abelian, resp. K3.
We refer to \cite{Sim} or to \cite{HL} for the definition
and the basic properties of the Simpson \mbox{(semi-)}stable sheaves.
Let $M_Y^{H,s}(v)$ (resp. $M_Y^{H,ss}(v)$) denote
the moduli space of Simpson stable (resp. semistable)
sheaves $\FFF$
on $Y$ with respect to an ample class $H$ with
Mukai vector $v(\FFF)=v$.
According to Mukai (\cite{Mu-1}, \cite{Mu-2}, see also \cite{HL}),
$M_Y^{H,s}(v)$, if nonempty,
is smooth of dimension $v^2 +2$ and carries
a holomorphic symplectic structure. 

\begin{theorem}\label{many-authors}
If $Y$ is a $K3$ surface, then a nonempty moduli space
$M_Y^{H,s}(v)$ is an irreducible symplectic variety whenever
it is compact or, equivalently, projective.  Moreover, $M_Y^{H,s}(v)$
is compact if $v$ is primitive and
$H$ is a sufficiently generic ample class in the
K\"ahler cone of $Y$. The last condition means that $H$ does
not lie on a certain discrete family of walls in the K\"ahler cone,
and it is automatically verified if $\Pic Y\simeq \ZZ$.
\end{theorem}

\begin{proof}
This summarizes the results of several papers: \cite{Mu-1}, \cite{Mu-2},
\cite{Hu-1}, \cite{Hu-2}, 
\cite{O'G-1} and \cite{Y-1}.
See also \cite{HL}, 6.2.5 and historical comments to Chapter 6.
\end{proof}

In
particular, the following statement holds:

\begin{corollary}\label{irred-symp-K3}
Let $Y$ be a K3 surface with $\Pic Y\simeq \ZZ$ and
$v$ a primitive Mukai vector. Assume that $M=M_Y^{H,s}(v)$ is nonempty. Then $M$
is an irreducible symplectic variety of dimension $v^2 +2$.
If $v$ is, moreover, isotropic, then $M$ is a K3 surface.
\end{corollary}

There are similar results for the case when
$Y$ is abelian \cite{Y-1}, however in this case not $M_Y^{H,s}(v)$ itself
is irreducible symplectic, but the fiber of its Albanese map
$M_Y^{H,s}(v)\rar Y\times\hat Y$. 
The papers cited above prove also that all the irreducible
symplectic varieties obtained in this way are deformation
equivalent to Beauville's examples.

Let $X$ be a K3 surface containing a nef curve $C$ of degree $2d-2$,
$d\geq 2$. 
Then the $d$-th punctual Hilbert scheme 
admits a dominant rational map $\theta:\Xdd\lrdash|C|\simeq\PP^{d}$ sending
$\xi\in\Xdd$ to the generically unique curve $C_\xi\in|C| $ containing $\xi$.
If $|C|$ embeds $X$ into $\PP^d$, then $\theta$ can be described as follows:
\begin{equation}\label{theta}
\theta :\Xdd\lrdash\PP^{d\dual}\ ,\ \ \ \xi \mapsto \langle \xi\rangle_{\PP^{d}}\ .
\end{equation}
Here $\langle \xi\rangle_{\PP^{d}}$
denotes the linear span of a subscheme $\xi\subset X$ in the embedding
into $\PP^{d}$, which is generically a hyperplane in $\PP^{d}$, that is
a point of the dual projective space $\PP^{d\dual}$.

For a $2d$-dimensional symplectic variety $V$,  we will call
a morphism $\pi:V\rar B$ a {\em Lagrangian fibration} if it is surjective and
its generic fiber is a connected Lagrangian subvariety of $V$,
that is a $d$-dimensional subvariety such that the restriction of
the symplectic form to it is zero. By the classical Liouville's theorem,
the generic fiber is then an abelian variety. 

\begin{theorem}[Matsushita, \cite{Mat-1}, \cite{Mat-2}]\label{matsushita-thm}
Let $V$ be an irreducible $2d$-dimensional symplectic variety, and $\pi:V\rar B$ a
surjective morphisme with connected fibers. Then $\dim B=d=\frac{1}{2}\dim V$
and $f$ is a Lagrangian fibration. If, moreover, $B$ is nonsingular,
then $R^if_*\OOO_V\simeq \Omega^i_B$ for all $i\geq 0$, and $B$ has
the Hodge numbers of a projective space.
\end{theorem}

We will call a rational
map $\pi:V\lrdash B$ a {\em rational Lagrangian fibration}, if there exists
a birational map $\mu :W\dasharrow V$ of $2d$-dimensional symplectic varieties
such that $\pi\mu$ is a regular Lagrangian fibration.
Such a $\pi$ is dominant and its generic fiber is
a connected Lagrangian subvariety of $V$, birational to an abelian variety.
The above map $\theta$ is a rational
Lagrangian fibration.
Its fibers are birational to symmetric powers $C_\xi^{(d)}$, and hence
to Jacobians of the genus-$d$ curve $C_\xi$. This birationality is globalized
as follows.

Let $J=J^dX$ be the relative compactified Jacobian of the linear
system $|C|$. It can be defined as the moduli space $M_X(0,[C],1)$
of Simpson-semistable torsion sheaves on $X$ with Mukai vector $(0,[C],1)$
\cite{Sim}, \cite{LeP-2}. To speak about semistable sheaves,
one has to fix a polarization $H$ on $X$, so it is better to write
$J^{d,H}X$ to explicitize the dependence on the polarization.
If all the curves in $|C|$ are reduced and irreducible, then
every semistable sheaf is stable, hence $J$ is smooth and its definition
does not depend on polarization. In this case $J$ can be
equally understood as the moduli space of simple
sheaves \cite{LeP-1}, \cite{Mu-1} with given Mukai vector. 
Let $\psi :J\lra |C|\simeq \PP^d$
be the natural map sending each sheaf to its support. It is known
that $J$ is an irreducible symplectic variety and $\psi$ is
a Lagrangian fibration \cite{Beau-3}.
The varieties $\Xdd$ and $J$ are related by a generalized Mukai flop
$\mu :\Xdd\lrdash J$ introduced by Markman in \cite{Mar}. It sends $\xi\in\Xdd$
to the same subscheme $\xi$ considered as a degree-$d$ divisor on the curve
$C_\xi$, and we have $\theta=\psi\circ\mu$.

According to Huybrechts (Lemma (2.6) of \cite{Hu-2}), any birational map
between irreducible symplectic varieties induces a Hodge isometry of their integer
second cohomology lattices  $H^2(\cdot,\ZZ)$
equiped with the Bogomolov--Beauville form. Thus we have the isomorphisms of the
Bogomolov--Beauville lattices $H^2(J,\ZZ)\simeq H^2(\Xdd,\ZZ)$
and $\Pic (J)\simeq \Pic(\Xdd)$.

\begin{lemma}\label{fminuse}
Let $X$ be a $K3$ surface with an effective divisor class $f_{2d-2}$
such that all the curves in the linear system $|f_{2d-2}|$ are reduced 
and irreducible. Let $D$ be a divisor on $\Xdd$ with class $f_{2d-2}-e$.
Then $h^0(\OOO (D))=d+1$ and  $\theta$ is given by the 
complete linear system $|D|$.
\end{lemma}

\begin{proof}
Let $\LLL=\psi^*\OOO_{\PP^{d}}(1)$. Then the self-intersection $(\LLL)^{2d}$
is 0, hence $q_J(\LLL)=0$, where $q_J$ denotes the Bogomolov--Beauville form
on $H^2(J,\ZZ)$. 
By the Riemann--Roch Theorem for
hyperk\"ahler manifolds \cite{Hu-1}, if $V$ is a
hyperk\"ahler manifold $V$ of dimension $2d$
with Bogomolov--Beauville form $q_V$, then $\chi (\OOO_V (D))=
\chi (\OOO_V)=d+1$ for any divisor $D$ such that $q_V(D)=0$.
Hence $\chi (\LLL)=d+1$.
By Matsushita's theorem \cite{Mat-2},
$R^i\psi_*\OOO_J\simeq \Omega^i_{\PP^d}$. Applying the Leray
spectral sequence and the Bott formula, we obtain for $\LLL=\psi^*(\OOO(1))$:
$h^0(\LLL)=d+1$,
$h^i(\LLL)=0$ for $i>0$.

The map $\theta=\psi\circ\mu$
is given by a subsystem of $|\mu^*\LLL|$. As $\mu^*$ is an isometry
of Bogomolov--Beauville lattices $(H^2(J,\ZZ), q_J)$ and $(H^2(\Xdd,\ZZ), q)$,
$q(\mu^*\LLL)=q_J(\LLL)=0$. For $X$ as in the hypothesis, 
it is possible that there are many primitive divisors $D$ with $q(D)=0$
on $\Xdd$. Let us consider a deformation of $X$ to a surface with
$\Pic X\simeq \ZZ$,
polarized by a divisor $f_{2d-2}$ of degree $2d-2$. Then 
the map $\theta$ deforms with $X$, $\rk\Pic \Xdd=2$ and the only two primitive
effective classes with square 0 on $\Xdd$ are $f_{2d-2}-e$,\ \ 
$f_{2d-2}+e$. The latter has negative intersection with the generic fiber
$\PP^1$ of the Hilbert--Chow map $\Xdd\lra X^{(d)}$, hence has the whole
exceptional divisor $E$ in its base locus,
which is not the case for $\theta$, so the linear system
defining $\theta$ is a subsystem of $|f_{2d-2}-e|$.

As $\mu$ is an isomorphism in codimension 1, $h^0(\mu^*\LLL)=h^0(\LLL)=d+1$,
and this ends the proof.
\end{proof}

As follows from Proposition \ref{H2-for-Xdd}, if $X$ is
a K3 surface with a curve class $f_{2d-2}$ of degree $2d-2$,
then the Bogomolov--Beauville
form on the Picard lattice of $\Xdd$ represents zero. The 
classes $\pm f_{2d-2}\pm e \in\Pic (\Xdd) $  are primitive with
square 0, and one of them,
namely $f_{2d-2}-e$, defines a rational Lagrangian fibration.
Thus Lemma \ref{fminuse} provides an example 
illustrating the following conjecture:

\begin{conjecture}\label{main-conj}
Let $X$ be a K3 surface. Then $\Xdd$ admits a
rational Lagrangian fibration with base $\PP^d$ if and only if
the Bogomolov--Beauville form of the Picard lattice of $\Xdd$ represents zero.
In this case, there exists $m\geq 1$ and an effective curve
class $f_{(2d-2)m^2}$ of degree $(2d-2)m^2$ on $X$ such that
the linear system defining the rational Lagrangian fibration
map on $\Xdd$ is of the form $f_{(2d-2)m^2}-me$.
\end{conjecture}

In Section \ref{fmtransform}, we will prove the conjecture
for K3 surfaces with Picard group $\ZZ$, generated by a curve class of
degree $(2d-2)m^2$,
and we will identify the Lagrangian fibration as
a torsor under the relative Jacobian of a linear
system of curves on some other K3 surface. The structure of the torsor is
defined by twisting by an element of the Brauer group of
the second K3 surface.

\section{Twisted sheaves and twisted Fourier-Mukai transforms}
\label{twisted}

The sheaves twisted by an element of the Brauer group
were introduced by C\u ald\u araru \cite{Cal-1}, \cite{Cal-2}.

The cohomological Brauer group $\Br'(Y)$ of a scheme $Y$ is defined
as $H^2_\et(Y,\OOO_Y^*)$. The Brauer group $\Br(Y)$
is the union of the images of
$H^1(Y,PGL(n))$ in $\Br'(Y)$ for all $n$. For a smooth curve,
$\Br(Y)=\Br'(Y)=0$. For a smooth surface, $\Br(Y)=\Br'(Y)$.
If $Y$ is a K3 surface, then $\Br(Y)\simeq
\Hom_\ZZ(T_Y,\QQ/\ZZ)$, where $T_Y$ is the transcendental
lattice of $Y$, defined as the orthogonal complement of
$\Pic Y$ in $H^2(Y,\ZZ)$.

Let $\alp\in\Br (Y)$ be represented by a \v Cech 2-cocycle
$(\alp_{ijk})_{i,j,k\in I}$, $\alp_{ijk}\in \Gamma(U_i\cap U_j\cap U_k, \OOO_Y^*)$
on some open covering  $\{U_i\}_{i\in I}$. An $\alp$-twisted
sheaf $\FFF$ on $Y$ is a pair
$(\{\FFF_i\}_{i\in I},\{\phi_{ij}\}_{i,j\in I})$,
where $\FFF_i$ is a sheaf on $U_i$ ($i\in I$) and
$$
\phi_{ij}:\FFF_j|_{U_i\cap U_j}\lra \FFF_i|_{U_i\cap U_j}\qquad (i,\ j\in I)
$$
are isomorphisms of sheaves with the following three properties:
$$
\phi_{ii}=\id, \ \ \phi_{ji}=\phi_{ij}^{-1}, \ \ 
\phi_{ij}\circ \phi_{jk}\circ \phi_{ki}=\alp_{ijk}\cdot\id .
$$

The $\alp$-twisted sheaves form an abelian category
$\Mod (Y,\alp)$. 
It has enough injectives and enough $\OOO_Y$-flat
sheaves.
An $\alp$-twisted sheaf is coherent if all the sheaves $\FFF_i$
are. The $\alp$-twisted coherent sheaves form an abelian
category $\Coh (Y,\alp)$. If $\FFF$ is an $\alp$-twisted sheaf,
and $\GGG$ an $\alp'$-twisted sheaf, then $\FFF\otimes \GGG$
is an $\alp\alp'$-twisted sheaf, and $\HOM(\FFF,\GGG)$ is
an $\alp^{-1}\alp'$-twisted sheaf. If $f:X\rar Y$ is a morphism,
then $f^*\FFF\in\Mod (X,f^*\alp)$, and for any 
$\GGG\in \Mod (X,f^*\alp)$, \ \ $f_*\GGG\in \Mod (Y,\alp)$.

We denote by $\DDB (Y,\alp)$ the derived category
$\DDB^b_\coh (\Mod (Y,\alp))$ of bounded complexes
of $\alp$-twisted sheaves with coherent cohomology.
By a standard machinery one defines the derived functors
$\RD f_*$, $\Lotimes$, $\LD f^*$.

To define the Chern character on $\Coh (Y,\alp)$,  C\u ald\u raru
fixes some
\mbox{$\alp^{-1}$-twisted} locally free sheaf $\EEE$ on $Y$. Then $\FFF\otimes\EEE$
is a usual (untwisted) sheaf, so one can define the modified Chern character
$
\ch_\EEE(\FFF):=\frac{1}{\rk\EEE}\ch(\FFF\otimes\EEE)
$
and the associated Mukai vector  $v_\EEE(\FFF)=\ch_\EEE(\FFF)\sqrt{\Td(Y)}$.
Remark that such an $\EEE$ exists only with rank which is a
multiple of the order of $\alp$ in the Brauer group.
These definitions depending on $\EEE$
are not the ones best suited for application to the Fourier--Mukai transform.
Huybrechts and Stellari \cite{HS} tensor $\FFF$
by an \mbox{$\alp^{-1}$-twisted}
$C^\infty$ line bundle in place of 
an \mbox{$\alp^{-1}$-twisted} holomorphic locally free sheaf $\EEE$ of higher rank.
Their construction depends on a so called $B$-field $B\in H^2(Y,\QQ)$.
The latter can be defined for any $Y$ for which $H^3(Y,\ZZ)=0$ as a lift of $\alp$
via the surjection in the exact triple
$$
0\rar \Pic Y\otimes\QQ/\ZZ\rar H^2(Y,\QQ/\ZZ)\rar H^2(Y,\OOO_Y^*)_\tor\rar 0
$$
composed with the natural map $H^2(Y,\QQ)\rar H^2(Y,\QQ/\ZZ)$. We will write
$\alp=e^{2\pi iB}$.

Let $L_B$ denote the $C^\infty$ line bundle on
$Y$ with transition functions $e^{2\pi i\beta_{ij}}$, where $(\beta_{ij})$
is a $C^\infty$ 1-cochain whose coboundary is some 2-cocycle
$B_{ijk}\in\Gamma(U_i\cap U_j\cap U_k, \QQ) $ representing $B$. Then
the twisted Chern character is $\ch_B(\FFF)=\ch (\FFF\otimes L_B^{-1})$
and the twisted Mukai vector is $v_B(\FFF)=\ch_B(\FFF)\sqrt{\Td(Y)}$.

Let now $X$ and $Y$ be smooth projective varieties, $\alp\in\Br(Y)$,
and $\PPP^\fatdot\in\DDB(X\times Y, \pi_Y^*\alp^{-1})$, where
$\pi_X$, $\pi_Y$ denote  the  projections of $X\times Y$
to the two factors. The twisted Fourier--Mukai transform
$$
\Phi_{Y\rar X}^{\PPP^{\textstyle\cdot}}:
\DDB(Y,\alp)\lra \DDB(X)
$$
is defined by
$$
\Phi_{Y\rar X}^{\PPP^{\textstyle\cdot}}(K^\fatdot)=
\pi_{X*}(\pi_Y^*(K^\fatdot)\Lotimes\PPP^\fatdot).
$$
It can be pushed down to the cohomology level in a natural way to give a
map $\phi:H^*(Y,\QQ)\rar H^*(X,\QQ)$ so that
the following diagram is commutative:
$$\begin{CD}
\DDB(Y,\alp)@>\Phi_{Y\rar X}^{\PPP^{\textstyle\cdot}}>> \DDB(X)\\
@V{v_B(\cdot)}VV @VV{v(\cdot)}V \\
H^*(Y,\QQ)@>\phi>> H^*(X,\QQ)
\end{CD}
$$
The Grothendieck--Riemann--Roch Theorem implies the following expression
for $\phi=\phi_{Y\rar X}^{\PPP^{\textstyle\cdot}\!,\:B}$:
$$
\phi_{Y\rar X}^{\PPP^{\textstyle\cdot}\!,\:B}(\gamma)=
\pi_{X*}\left(\pi_Y^*(\gamma)\cdot\ch_{\pi_Y^*B}(\PPP^\fatdot)\cdot\sqrt{\Td(X\times
    Y)}\right)\ .
$$
If $\alp=1$, $B=0$, we get the usual Fourier--Mukai
transform introduced by Mukai, and 
$\phi_{Y\rar X}^{\PPP^{\textstyle\cdot}\!,\:0}$ is denoted by
$\phi_{Y\rar X}^{\PPP^{\textstyle\cdot}}$.
Remark that the cohomological Fourier--Mukai
transform $\phi$ does not respect the grading of cohomology, neither the ring
structure given by the cup product.

Let now $X$ be a K3 surface and $v=(r,L,s)\in H^*(X,\ZZ)$ a  primitive
Mukai vector with algebraic $L\in H^2(X,\ZZ)$ such that \mbox{$v^2:=(v,v)=0$}.
Let $H$ be a sufficiently generic ample
class, so that $M=M_X^{H,s}(v)$ is compact. Then, by Corollary \ref{irred-symp-K3},
$M$ is a K3 surface. Denote $m=\GCD(r,L\cdot \gamma ,s)_{\gamma\in \Pic X}$. 
If $m=1$, then, according to the Appendix to  \cite{Mu-2},
$M$ is a fine moduli space, that is, there exists a universal sheaf
$\PPP$ on $X\times M$. It is defined by the condition that for any $t\in M$,
the isomorphism class of stable sheaves corresponding to $t$ is represented by
the restriction $\PPP_t=\PPP|_{X\times t}$. The same Mukai's argument
shows that if $m>1$, then there exists an element $\alp\in\Br(M)$ of
order $m$ and  $\pi_M^*\alp^{-1}$-twisted universal
sheaf $\PPP$ on $X\times M$.
We will not consider separately the case  $m=1$. It is a particular
case of the next theorem corresponding to $\alp=1$, $B=0$, though
historically, it was proved by Mukai in \cite{Mu-2}, \cite{Mu-3}
before C\u ald\u araru's work on twisted sheaves. To state the theorem,
we need to introduce a new weight-2 integral Hodge structure on the total
cohomology $H^*(M,\ZZ)$ of $M$. This Hodge structure is determined
by a $B$-field lifting $\alp$, which we will
fix once and forever, and will be denoted by
$\tilde H_B(M)$. The integer structure is defined by
$$
\tilde H_B(M,\ZZ)=(\exp B\cdot H^*(M,\QQ))\cap H^*(M,\ZZ), \ \exp B:=
1+B+\frac{B^2}{2},
$$
and the Hodge decomposition over $\CC$ by
$$
\tilde H_B^{1,1}(M)=\exp B\cdot (H^0(M)\oplus
H^{1,1}(M)\oplus H^4(M)),
$$
$$
\tilde H_B^{2,0}(M)=\exp B\cdot H^{2,0}(M),\ \ \tilde H_B^{0,2}(M)=\exp B\cdot H^{0,2}(M).
$$
The lattice $\tilde H_B(M,\ZZ)$ is endowed with Mukai's form
(\ref{mukai-form}).
We will also consider the total cohomology of $X$ with a
weight-2 Hodge structure $\tilde H(X,\ZZ)$ constructed in the same way, but with $B$-field equal to zero.

The following theorem holds.

\begin{theorem}[Mukai, C\u ald\u araru, Huybrechts--Stellari]\label{MCHS}
Under the hypotheses and with the notation of the previous paragraph, \sloppy\ 
let $\PPP$ be a \mbox{$\pi_M^*\alp^{-1}$-twisted}  universal sheaf on $X\times M$ and
$\PPP^\dual:=\R\HOM (\PPP,\OOO_{X\times M})$.
The following assertions are verified:

(i)\ \mbox{$
\Phi_{M\rar X}^{\PPP}:
\DDB(M,\alp)\lra \DDB(X)
$}\sloppy\ 
is an equivalence of categories with inverse
\mbox{$
\Phi_{X\rar M}^{\PPP\sdual}:
\DDB(X)\lra \DDB(M,\alp)
$}.

(ii) $\phi=\phi_{X\rar M}^{\PPP{\sdual}\!,\:B}$ is defined over $\ZZ$ and is
a Hodge isometry 
$\tilde H(X,\ZZ)\isoto\tilde H_B(M,\ZZ)$.

(iii) We have $\phi (v)=(0,0,1)$, $\phi (v^\perp)\subset (0,*,*)$, and
the induced map
$$
\bar\phi:v^\perp/v\rar H^2(M,\ZZ),\ \ w\mapsto [\phi (w)]_{H^2(M)}
$$
is a Hodge isometry. Here $[\cdot]_{H^2(M)}$ denotes the
$H^2(M)$-component of an element of $H^*(M)$, and the orthogonal
complement $v^\perp$ is taken in
$\tilde H(X,\ZZ)$.
\end{theorem}

Yoshioka in \cite{Y-2} defined the notion of (semi)-stability of $\alp$-twisted
sheaves and constructed the moduli spaces of $\alp$-twisted sheaves. Let
$Y$ be a K3 surface polarized by an ample class $H$, $B\in H^2(Y,\QQ)$
and $\alp =e^{2\pi iB}\in \Br(Y)$. Let $v\in \tilde H_B^{1,1}(M,\ZZ)$.
We will denote by $M_{Y,B}^{H,s}(v)$ the moduli space
of stable $\alp$-twisted sheaves $\FFF$ on $Y$ with $v_B(\FFF)=v$.
The next theorem is a generalization of
Theorem \ref{many-authors} to twisted sheaves.

\begin{theorem}[Yoshioka]\label{yoshioka-thm}
In the hypotheses of the previous paragraph, assume that
$H$ is sufficiently generic, that $v$ is primitive, and that
$M=M_{Y,B}^{H,s}(v)$ is nonempty. Then $M$ is an irreducible
symplectic manifold of dimension $v^2+2$.
\end{theorem}

\section{The K3 moduli space $M(m,H,(d-1)m)$}\label{modspace}

Let $X$ be a K3 surface such that $\Pic X=\ZZ\cdot H$, where
$H$ is ample and $H^2=(2d-2)m^2$ for some
$m\geq 2$, $d\geq 2$. Denote by  $\eta_X$
the positive generator
of $H^4(X,\ZZ)$. Let $M=M_X(m,H,(d-1)m)$ be
the moduli space of semistable sheaves on $X$ with
Mukai vector $v=(m,H,(d-1)m)=m+H+(d-1)m\eta_X$.

\begin{lemma}
The following statements hold:

(i) $M$ is a K3 surface with $\Pic M\simeq \ZZ$.

(ii) Every
semistable sheaf with Mukai vector $v$ is $\mu$-stable and locally free.

\end{lemma}

\begin{proof}
As $v$ is primitive and $v^2=0$, (i) follows from Theorem \ref{MCHS}, (iii).
Further, $\Pic X=\ZZ H$, hence every
semistable sheaf $\EEE$ with $c_1(\EEE)=H$ is $\mu$-stable.
The local freeness of $\EEE$ follows
from  \cite{Mu-2}, \mbox{Proposition 3.16.}
\end{proof}

Now we will describe Serre's construction, which permits to
obtain all the vector bundles from
$M$ as some sheaf extensions.
For a \mbox{$0$-dimensional} subscheme $Z\subset X$, the
number $\delta (Z)=h^1(\III_Z(1))$ is called the index of speciality
(or $\OOO (1)$-speciality)
of $Z$; it is equal to the number of independent
linear relations between the points of $Z$.
In a more formal way, $\delta (Z) =
l(Z)-1-\dim \langle Z\rangle$, where
$l(Z)$ stands for the length of $Z$. Following Tyurin \cite{Tyu-1}, \cite{Tyu-2},
we will call $Z$ stable if $\delta (Z')<\delta (Z)$
for any $Z'\subset Z$, $Z'\neq Z$. 

\begin{lemma}
Let $Z\subset X$ be a stable $0$-dimensional subscheme 
of degree $c=(d-1)(m^2-m)+m$ with  $\delta (Z)=m-1$.
Define a sheaf $\EEE =\EEE_Z$ as the middle term of
the exact triple
\begin{equation}\label{reider}
0\lra H^1(\III_Z(1))\otimes\OOO_X\xrightarrow{\alp}
\EEE\xrightarrow{\beta}\III_Z(1)\lra 0
\end{equation}
whose extension class is the identity map
on $H^1(\III_Z(1))$ under the canonical isomorphism
$\Ext^1(\III_Z(1),H^1(\III_Z(1))\otimes\OOO_X)=
\End(H^1(\III_Z(1))$.
Then $\EEE$ is a stable locally free sheaf with Mukai vector $v=\mbox{$(m,H,(d-1)m)$}$
and $h^0(\EEE)=dm$, $h^1(\EEE)=h^2(\EEE)=0$.
\end{lemma}

\begin{proof}
The local freeness follows from \cite{Tyu-2}, Lemma 1.2.
The assertions on the cohomology and the Mukai vector of $\EEE$
are obvious.
As $\Pic (X)=\ZZ H$ and $c_1(\EEE)=H$, it is enough
to prove that $\EEE$ is semistable.

Assume that $\EEE$ is unstable and $\MMM$ is the maximal
destabilizing subsheaf of $\EEE$. Then $c_1(\MMM )=nH$
with $n\geq 1$. Let $i$ be the inclusion $\MMM\lra\EEE$.
If $\beta\circ i=0$, then $\alp^{-1}$ maps $\MMM$ into the
trivial vector bundle $H^1(\III_Z(1))\otimes\OOO_X\simeq \OOO_X^{m-1}$,
which is impossible. Hence $\beta\circ i(\MMM )$
is a rank-1 subsheaf of $\III_Z(1)$, which we can
represent as $\III_W(1)$ for some subscheme $W\subset X$
containing $Z$. Thus we
have the exact triple
$$
0\lra\MMM'\lra\MMM\lra\III_W(1)\lra 0.
$$
If we assume that $W$ has one-dimensional
components, then $c_1(\III_W(1))\leq 0< \mu (\MMM )$,
which contradicts the semistability of $\MMM$. Hence
$W$ is 0-dimensional.
We have $\MMM'\subset\OOO_X^{m-1}\subset\EEE$, so
$c_1(\MMM')=n-1\leq 0$ and $n=1$.
As $\MMM$ is maximal, it is saturated in $\EEE$,
hence so is $\MMM'$, and then  $\MMM'\simeq\OOO_X^{k-1}$
is a trivial subbundle of $\OOO_X^{m-1}$. We obtain the
following commutative diagram with exact rows and columns:
\begin{equation}\label{CD1}
\begin{CD}
 @. 0 @. 0 @. 0 @. \\
@. @VVV @VVV @VVV @. \\
0 @>>>\OOO_X^{k-1}  @>>> \MMM @>>>
\III_W(1) @>>> 0 \\
@. @VVV @VVV @VVV @. \\
0 @>>>\OOO_X^{m-1}  @>>> \EEE @>>> \III_Z(1) @>>> 0 \\
@. @VVV @VVV @VVV @. \\
0 @>>> \OOO_X^{m-k} @>>> \NNN @>>> \III_{W,Z} @>>> 0 \\
@. @VVV @VVV @VVV @. \\
 @. 0 @. 0 @. 0 @. \\
\end{CD}
\end{equation}
For any Artinian $\OOO_X$-module $M$, we have
$\Ext^1(M, \OOO_X)\simeq\Ext^1(\OOO_X,M)^\dual=H^1(X,M)^\dual=0$.
Applying this to $M=\III_{W,Z}$, we see that
$\Ext^1(\III_{W,Z},\OOO_X^{m-k})=0$
and $\NNN$ has torsion whenever $\III_{W,Z}\neq 0$, which contradicts
the fact that $\MMM$ is a saturated subsheaf of $\EEE$. Hence
 $\III_{W,Z}= 0$ and $Z=W$.

We have the following monomorphism of extensions of sheaves:
$$
\begin{CD}
0 @>>>\OOO_X^{k-1}  @>>> \MMM @>>>
\III_Z(1) @>>> 0 \\
@. @VVV @VVV @| @. \\
0 @>>>\OOO_X^{m-1}  @>>> \EEE @>>> \III_Z(1) @>>> 0 \\
\end{CD}
$$
The next lemma implies that $k=m$, hence
the monomorphism is an isomorphism,
$\MMM=\EEE$, and $\EEE$ is semistable.
\end{proof}

\begin{lemma}
Let $V$ be a smooth variety, $F_1,F_2,G$ sheaves on $V$,
$$
0\lra F_i\lra \EEE_i\lra G\lra 0\ ,\qquad i=1,\ 2,
$$
two sheaf extensions with classes $e_i\in
\Ext^1(G,F_i)$. Then a homomorphism of sheaves $\phi : F_1\lra F_2$
extends to a morphism of extensions
$$
\begin{CD}
0 @>>>F_1  @>>> E_1 @>>>
G @>>> 0 \\
@. @V{\phi}VV @VVV @| @. \\
0 @>>>F_2  @>>> E_2 @>>> G @>>> 0 \\
\end{CD}
$$
if and only if $e_2=\phi\circ e_1$.
\end{lemma}

\begin{proof}
Standard; compare to \cite{Mac}, Proposition III.1.8.
\end{proof}

The vector bundle $\EEE=\EEE_Z$ defined by (\ref{reider})
is the result of {\em Serre's construction} applied to the
0-dimensional susbscheme $Z$.

\begin{remark}
To define  $\EEE_Z$, one can replace the identity
extension class by any linear automorphism  $e$ of $ H^1(\III_Z(1))$.
By the lemma, this will provide an isomorphic sheaf  $\EEE_Z$.
If $e$ is not of maximal rank, then $\EEE_Z$
is not locally free. By \cite{Tyu-2}, \S\ 1, $\EEE_Z$
will also be non-locally-free if $Z$ is not stable.
\end{remark}

The description of properties of $\EEE$ in terms of $Z$ is particularly
simple if $Z$ is contained in a smooth hyperplane section
$C\in |H|$ of $X$, in which case we can consider $Z$
as a divisor on $C$. By the geometric Riemann--Roch
Theorem, the index of speciality of $Z$ according to our definition
is the same as the index of speciality of the divisor $Z$ on $C$.
Thus $\delta (Z)=m-1$ if and only if $Z$ belongs to
a linear series $g_c^{m-1}$ on $C$. We denote by
$W_c^r(C)$ the Brill--Noether locus
of linear series $g_c^r$ on $C$ and by $G_c^r(C)$ the union
the corresponding linear series as a subvariety of $\Div^c(C)$.

\begin{lemma}\label{fitting}
Let $C\in |H|$ be a smooth curve and $Z\in\Div^c(C)$ such that
$\delta (Z)=m-1$. Let $\EEE_Z$ be defined by the extension
(\ref{reider}). Then the following assertions hold:

(i) $\EEE_Z$ is locally free if and only if the linear system
$|Z|$ is base point free.

(ii) $\EEE_Z$ is globally generated if and only if the linear system
$|K_C-Z|$ is base point free.

(iii) $\EEE_Z$ fits into the exact triple
$$
0\lra \OOO_X^m\lra\EEE_Z\lra\OOO_C(K-Z)\lra 0\ ,
$$
where $K=K_C$ is the canonical class of $C$,
and $\EEE_Z\simeq\EEE_{Z'}$ for any $Z'\in |Z|$.
\end{lemma}

\begin{proof}
See \cite{Tyu-1}, Lemma 3.4 or
\cite{Mor}, Section 5.
\end{proof}

For any $t\in M$ we will denote by $\EEE_t$ the rank-$m$ sheaf
on $X$ represented by $t$.

\begin{proposition}\label{globgen}
Let $t\in M$ be generic and $\EEE=\EEE_t$.
The following properties are verified:

(i) $\EEE$ is obtained from a divisor $Z\in G_c^{m-1}(C)$
for some smooth hyperplane section $C$ of $X$ by the
construction of Lemma {\em \ref{fitting}}.

(ii) $\EEE$ is globally generated.
\end{proposition}

\begin{proof}
As $X$ is generic, any of its smooth hyperplane sections $C$
is Brill--Noether generic by \cite{L}. As the Brill--Noether
number $\rho_{c}^{m-1}=1$, there exists a divisor
$Z\in g_{c}^{m-1}$ on $C$. Let $\EEE =\EEE_Z$ be the corresponding
vector bundle.
As $\rho_{c-1}^{m-1}<0$ and $\rho_{c+1}^{m}<0$, the linear
systems $|Z|$ and $|K-Z|$ are base point free, so $\EEE =\EEE_t$
for some $t\in M$ and $\EEE$ is globally generated by Lemma \ref{fitting}.
The conditions (i), (ii) are open, so the proposition is proved.
\end{proof}


For a future use, we will prove the following lemma.

\begin{lemma}\label{technical}
Let $t\in M$ be generic and $\EEE=\EEE_t$.
The following properties are verified:

(i) Let $s_1,\ldots,s_m$ be $m$ linearly independent global
sections of $\EEE$. Then $s_1,\ldots,s_m$ generate $\EEE$
at generic point of $X$.

(ii) The vanishing at a generic point $z\in X$ imposes precisely $m=\rk \EEE$
independent linear conditions on a section of $\EEE$. For any
$k=1,\ldots,d$ and for $k$ generic points $z_1,\ldots ,z_k$ of $X$,
we have $h^0(X,\EEE\otimes \III_\xi )=m(d-k)$, where
$\xi=z_1+\ldots+z_k$.

(iii) Let $z_1,\ldots ,z_{d-1}$ be generic. Then any nonzero section $s$
of $\EEE$ vanishing at $z_1,\ldots ,z_{d-1}$ has exactly
$z_1+\ldots +z_{d-1}$ as its scheme of zeros.
\end{lemma}

\begin{proof}
(i) Let $\FFF$ be the saturate in $\EEE$
of the subsheaf generated by $s_1,\ldots,s_m$. As $\EEE$ is locally free, $\FFF$ is reflexive.
Any reflexive sheaf on a smooth surface is locally free, so $\FFF$ is locally free.
By the stability of $\EEE$, $\FFF$ has no subsheaves $\FFF'$ with
$c_1(\FFF')>0$. If $k=\rk\FFF<m$, then $s_1,\ldots,s_m$ are linearly
dependent at generic point of $\EEE$. Let us renumber the $s_i$
in such a way that $s_1,\ldots,s_k$ are linearly independent
over $\CC (X)$. Then $s_1,\ldots,s_k$ define an inclusion $\OOO_X^{\oplus k}\into\FFF$. 
An inclusion of locally free sheaves of the same rank is either isomorphic,
or has a cokernel supported on a nonempty divisor. The second case is impossible,
as $c_1(\FFF)\leq 0$. Hence $\FFF\simeq \OOO_X^{\oplus k}$, and $s_i$
for $i=k+1,\ldots, m$ are linear combinations of $s_1,\ldots,s_k$ with constant
coefficients, which is absurd. Hence $k=m$ and $\FFF=\EEE$.

(ii) If $z_1$ is generic, then $H^0(\III_{z_1}\otimes\EEE)$ has
codimension $m$ in $H^0(\EEE)$ by (i). Choose
$s_1,\ldots,s_m\in H^0(\III_{z_1}\otimes\EEE)$ linearly independent.
Then for $z_2\in X$ generic, $s_1,\ldots,s_m$
span $\EEE_{z_2}$ by (i). Thus $z_1,z_2$ impose $2m$ conditions
on sections of $\EEE$. Iterating one gets (ii).

(iii) Choose $d-2$ generic points $z_1,\ldots,z_{d-2}$ of $X$ and a basis
of $2m$ sections $e_1,\ldots,e_m,s_1,\ldots,s_m$ vanishing at
$ z_1,\ldots,z_{d-2}$, as in (ii) with $k=d-2$.
Consider these sections on a sufficiently small open 
set $U$ on which both $e_1,\ldots,e_m$ and
$s_1,\ldots,s_m$ are bases of $\EEE$. Let $A=A(x)$ be a $m$ by $m$ matrix
of rational functions in $x\in X$ such that $s_i=Ae_i$, that is,
$s_i=\sum_j a_{ji}e_j$, where $A=(a_{ji})$ ($i,j=1,\ldots,m$).
Then $A$ is regular and nondegenerate for $x\in U$.
For any $y\in U$, the $m$ sections $\sigma_i=s_i-A(y)e_i$ 
form a basis of all the
sections of $\EEE$ vanishing at $ z_1,\ldots,z_{d-2}$ and $y$.
We may expand them in the basis $(e_j)$ with coefficients in $\CC (X)$:
$\sigma_i =\sigma_i (x)= (A(x)-A(y))e_i$. Here $x$ is a variable point of $U$,
so that $\sigma_i (x)$ denotes the value of $\sigma_i$ in the fiber
$\EEE_x$. Thus, any section $\sigma$ of $\EEE$ vanishing at $y$ 
can be written over $U$ in the form
$\sigma (x)=(A(x)-A(y))\vv$, where $\vv=v_1e_1+\ldots +v_me_m$,
$(v_1,\ldots,v_m)\in\CC^m$.
This $\sigma$ depends on $\vv$, $y$ as parameters, and we will denote it by $\sigma_{\vv,y}$.
By the stability of $\EEE$, no nonzero section of $\EEE$ vanishes along a curve,
so $y$ is an isolated zero of $\sigma_{\vv,y} $. Let us now fix $\vv\neq 0$ and let vary
both $x$ and $y$. The zeros of $\sigma_{\vv,y} $ are the solutions of the equation
$A(x)\vv=A(y)\vv$. The fact that $y$ is an isolated zero of $\sigma_{\vv,y}$
for any $y\in U$
implies that the fibers of the map $f:U\rar \CC^m$,
$y\mapsto A(y)\vv$ are 0-dimensional (here we interprete $\vv$
as the vector $(v_1,\ldots,v_m)\in\CC^m$). By Chevalley's Theorem,
$f(U)$ is a 2-dimensional constructible subset of $\CC^m$
in the Zariski topology.
By Bertini-Sard Theorem, $f(U)$ contains an open subset of noncritical values
of $f$. Thus there is a smaller Zariski open subset $U_0\subset U$ such that
$f|_{U_0}$ is locally holomorphically invertible (``locally'' in the classical
topology). Hence for any $y_0\in U_0$, the equation $f(x)=f(y_0)$ is locally
analytically equivalent to $x=y_0$, thus $y_0$ is a simple zero of $f(x)-f(y_0)$
and a simple zero of $\sigma_{\vv,y_0}$. 

We have proved that any nonzero section of $\EEE$ vanishing at $d-1$ generic points
$z_1,\ldots,z_{d-2},z_{d-1}=y_0$ of $X$ has a simple zero at $z_{d-1}$. By the
symmetry of the roles of the points $z_i$, all the $z_i$ are simple zeros.

\end{proof}

\section{Lagrangian fibration via Fourier--Mukai transform}
\label{fmtransform}

As in Section \ref{modspace},
let $X$ be a generic K3 surface of degree $(2d-2)m^2$.
Denote by $H$, resp. $\eta_X$
the positive generator
of  $\Pic X$, resp. $H^4(X,\ZZ)$. Let $M=M_X(m,H,(d-1)m)$ be
the moduli space of semistable sheaves on $X$ with
Mukai vector $v=(m,H,(d-1)m)$. Denote by $\eta_M$ the positive generator
of $H^4(M,\ZZ)$.

\begin{lemma}\label{generator}
Let $\PPP$ be a $\pi_M^*(\alp)^{-1}$-twisted universal sheaf
on \mbox{$X\times M$}\sloppy\ for some
$\alp\in\Br (M)$ of order $m$. Let $B\in H^2(M,\QQ)$ be such that $\alp=e^{2\pi iB}$ and
$\phi=\phi_{X\rar M}^{\PPP{\sdual}\!,\:B}$ the map
used in Theorem \ref{MCHS}.
Then $\Pic M\simeq\ZZ$, and for the positive generator $\hat H$ of
$\Pic M$ we have $$\hat H=\pm [\phi (1+(1-d)\eta_X)]_{H^2(M)},$$
where $[\cdot]_{H^2(M)}$ denotes the $H^2$-component of a cohomology class.
\end{lemma}

\begin{proof}
The intersection of $v^\perp\subset \tilde H(X,\ZZ)$ with $\tilde H^{1,1}(X)=
\CC\cdot 1+\mbox{$\CC\cdot H$}+\CC\cdot\eta_X$ is the lattice of rank 2 generated by
$v$ and $1+(1-d)\eta_X$. Hence, by Theorem \ref{MCHS}, the class of the hyperplane
section $\hat H$ of $M$ is equal to $\pm [\phi (1+(1-d)\eta_X)]_{H^2(M)}$.
\end{proof}

Let $\xi$ be a subscheme of length $d$ in $X$
and $\III_\xi\subset\OOO_X$ its ideal sheaf. 
Define
$$
C_\xi=\{ t\in M\mid h^0(X,\EEE_t\otimes \III_\xi )\neq 0\} ,
$$
where $\EEE_t$ denotes a rank-$m$ vector bundle whose isomorphism
class is represented by $t$.

\begin{proposition}\label{cxi}
The sign in the formula of Lemma \ref{generator} for $\hat H$ is plus, and
for generic $\xi\in X^{[d]}$,  $C_\xi$ is a curve
from the linear system $|\hat H|$.
\end{proposition}

\begin{proof}
Let $\Phi=\Phi_{X\rar M}^{\PPP\sdual}$ be the Fourier--Mukai transform associated to $\PPP$.
The Mukai vector \mbox{$1+(1-d)\eta_X$} is realized by either one
of the objects $\III_\xi$ or $\III_\xi^\dual =\mbox{$\R\HOM (\III_\xi,\OOO_X)$}$.
Hence 
$\phi (1+(1-d)\eta_X)={\ch_B}(\Phi(\III_\xi^\dual))\sqrt{\Td(M)}$,
where $\Phi(\III_\xi^\dual)=
\R \pi_{M*}(\pi_X^*\III_\xi^\dual
\Lotimes\PPP^\dual) $. As $\PPP$ is locally free, we may replace
$\Lotimes$ by $\otimes$. By the relative duality for $\pi_M$
(see \cite{Har-1}, p. 210 and \cite{Cal-2}, 2.7), and because the canonical sheaf
of $M$ is trivial, $\Phi(\III_\xi^\dual)[2]\simeq
(\R \pi_{M*}(\pi_X^*\III_\xi\otimes\PPP))^\dual$. The
Cohomology and Base Change Theorem (\cite{Har-2}, Theorem III.12.11)
reduces the computation of the latter to that of the
cohomology groups $H^i(\GGG_t)$, where $t\in M$,
$\GGG=\pi_X^*\III_\xi\otimes\PPP$, and $\GGG_t:=\GGG|_{X\times t}\simeq
\PPP_t\otimes\III_\xi$.
The cohomology of the sheaves $\PPP_t\otimes\III_\xi$ can be determined
from the exact triples
$$
0\lra H^1(\III_Z(1))\otimes\III_\xi\lra
\PPP_t\otimes\III_\xi\lra\III_{Z\cup\xi}(1)\lra 0,
$$
where we can assume that $\Supp Z\cap\Supp\xi=\emptyset$, and
$$
0\lra \III_{Z\cup\xi}(1)\lra \OOO_X(1)\lra \OOO_{Z\cup\xi}(1)\lra 0.
$$
By Serre duality,
$H^2(\PPP_t\otimes\III_\xi)=\Hom (\III_\xi,\PPP_t^\dual)^\dual=
H^0(\PPP_t^\dual)^\dual=0$. Hence $R^2\pi_{M*}\GGG=0$, and
$h^0(\PPP_t\otimes\III_\xi)=
h^1(\PPP_t\otimes\III_\xi)$, so both of them are different from zero
if and only if $t\in C_\xi $. 

Let us verify that 
$
C_\xi
$
is a proper closed subset of $M$ for generic $\xi$. 
It suffices to show that $h^0(X,\PPP_t\otimes \III_\xi )=0$
for generic $\xi\in X^{[d]}$ and generic $t\in M$.
Choose a generic $\EEE=\PPP_t$ and
$d$ generic points $z_1,\ldots ,z_d$ of $X$.
By Lemma \ref{technical}, (ii), $h^0(X,\EEE\otimes \III_\xi )=0$
for $\xi=z_1+\ldots+z_d$, which implies the result.

Thus, $C_\xi$ is a union of finitely many curves and isolated points
for generic $\xi$.
By Proposition 2.26 of \cite{Mu-2}, $C_\xi$ is of pure dimension~1 and
$R^0\pi_{M*}\GGG=0$.
As $R^i\pi_{M*}\GGG=0$ is nonzero only in odd degree $i=1$ and
$\EXT^p(R^1\pi_{M*}\GGG,\OOO_M)$ is nonzero only for odd $p=1$,
we have
$$
[\phi (1+(1-d)\eta_X)]_{H^2(M)}=
c_1(R^1\pi_{M*}\GGG)=h^0(\PPP_t\otimes\III_\xi)[C_\xi]
$$
for generic $t\in C_\xi$.
But by Lemma \ref{generator}, $\phi (1+(1-d)\eta_X)]_{H^2(M)}=\pm
\hat{H}$, hence $C_\xi\in |\hat{H}|$, $h^0(\PPP_t\otimes\III_\xi)$ is generically 1,
and the sign is plus.
\end{proof}

\begin{theorem} \label{main}
Let $X$ be a $K3$ surface with $\Pic X\simeq \ZZ$, and $H$ the ample generator
of $\Pic X$. Let $M=M_X^{H,s}(m,H,(d-1)m)$ be
the moduli space of $H$-stable sheaves on $X$ with
Mukai vector $v=(m,H,(d-1)m)$ ($d\geq 2, m\geq 2$) and $\hat H$
the ample generator of $\Pic M$.
Let $\PPP$ be a $\pi_M^*(\alp)^{-1}$-twisted universal sheaf
on \mbox{$X\times M$}\sloppy\ for some
$\alp\in\Br (M)$ of order $m$ and $B$ a lifting of $\alp$ in $H^2(M,\QQ)$.
Let 
$\Phi=\Phi_{X\rar M}^{\PPP{\sdual}}$ be the associated Fourier--Mukai transform, and
$\phi=\phi_{X\rar M}^{\PPP{\sdual}\!,\:B}$ its cohomological descent.
Denote by $w$
the Mukai vector $(1,0,1-d)$ of the sheaves $\III_\xi$ for $\xi\in X^{[d]}$,
so that $X^{[d]}=M_X^{H,s}(w)$. Then the following assertions hold:

(i) For generic $\xi\in  X^{[d]}$, the only nonzero cohomology of the complex
$\Phi (\III_\xi^\dual)$ is $h^2$, and $h^2\Phi (\III_\xi^\dual)$ is a
rank-1 torsion-free $\alp|_{C_\xi}$-twisted
sheaf on a curve $C_\xi\in |\hat H|$.

(ii) $\phi (w)=(0,\hat H, k)$ for some $k\in\ZZ$, and the moduli space
$V=M_{M,B}^{\hat H,s}(0,\hat H, k)$ is an irreducible symplectic manifold
of dimension $2d$. There is a birational
isomorphism $\mu:X^{[d]}\lrdash V$ defined by
$\xi\mapsto [h^2\Phi (\III_\xi^\dual)]$, 
where the brackets denote the isomorphism class of a sheaf.

(iii) The support of any sheaf $\LLL_t$ on $M$ represented by a point
$t\in V$ is a curve from the linear system $|\hat H|$, and the map
$f:V\rar |\hat H|\simeq \PP^d$, $t\mapsto \Supp \LLL_t$, is a Lagrangian fibration.
If we denote by $\{C\}$ the point of the projective space $\PP^d\simeq |\hat H|$
representing a curve $C$ from the linear system $|\hat H|$, then
the fiber $f^{-1}(\{C\})$ for generic $\{C\}\in |\hat H|$ is
isomorphic to the Jacobian of $C$.
\end{theorem}

\begin{proof}
(i) was verified in the proof of Proposition \ref{cxi}.

(ii) The equality $\phi (w)=(0,\hat H, k)$ follows from Proposition \ref{cxi} and
Theorem \ref{MCHS} (iii). The fact that $V$ is irreducble symplectic
will follow from Yoshioka's Theorem \ref{yoshioka-thm} as soon as we see that
$V$ is nonempty. But we have constructed stable sheaves represented by points of $V$
in part (i). Indeed, as $\Pic M=\ZZ\hat H$, any curve $C_\xi$ is irreducible,
and a rank-1 torsion
free (twisted or usual) sheaf on an irreducible curve is stable with respect 
to any polarization.

To prove the birationality of $\mu$, remark that $\Phi$ and the duality functor $D$
are equivalences of categories, so the composite functor $\Phi\circ D$ 
transforms non-isomorphic sheaves $\III_\xi$
into the complexes $\Phi (\III_\xi^\dual)$ that are non-isomorphic in $\DDB (M,\alp)$.
For generic $\xi$, the cohomology of the complex $\Phi (\III_\xi^\dual)$ is concentrated in degree 2,
hence the complex is quasi-isomorphic to $h^2\Phi (\III_\xi^\dual)[2]$.
Thus for generic $\xi\neq \xi'$, the sheaves $h^2\Phi (\III_\xi^\dual)$, $h^2\Phi (\III_{\xi'}^\dual)$
are non-isomorphic. This implies that $\mu$ is a generically injective rational map between
irreducible varieties of the same dimension, hence it is birational.

(iii) The fact that $f$ is a Lagrangian fibration is an obvious consequence of the above and of Matsushita's Theorem
\ref{matsushita-thm}. 

As $\Br(C)=0$ for a smooth curve $C$, $\alp|_{C}=0$ for any smooth $C\in |\hat H|$.
Hence the fiber of the support map over $C$ is isomorphic to the Jacobian $J(C)$.
This isomorphism is not canonical, for two different \v Cech 1-cochains $\beta$
such that $\check d(\beta)=\alp|_{C}$ may differ by a \v Cech 1-cocycle defining an
invertible sheaf $\LLL$ on $C$, and the corresponding
isomorphisms of $f^{-1}(\{C\})$ with $J(C)$ will differ
by a translation by the class of $\LLL$ in $J(C)$. Hence
$V$  represents a birational torsor
(biregular over the smooth curves $C\in |\hat H|$) under the relative Jacobian
$J$ of the linear system $|\hat H|$, and the generic fiber of $f$ is isomorphic to
$J(C)$ with $C\in |\hat H|$.
\end{proof}

\begin{corollary} \label{hminusme}
Let $X$ be a generic $K3$ surface of degree
$(2d-2)m^2$, $H$ the positive generator of $\Pic X$, and $h$
the divisor class in $\Pic (\Xdd)$
corresponding to $H$ under the isomorphism of Proposition \ref{H2-for-Xdd}.
Consider the rational map
$$
\pi :\Xdd\lrdash |\hat{H}|\simeq\PP^d,\ \ \xi\mapsto C_\xi\ .
$$
Then 
$\pi$ is defined by the complete linear system $|h-me|$ and
is a rational Lagrangian fibration. The fiber $\pi^{-1}(\{C\})$ for generic $C
\in|\hat{H}| $ is
birational to the Jacobian of $C$.
\end{corollary}

\begin{proof} In the notation of Theorem \ref{main}, $\pi=f\circ\mu$,
where $f$ is a Lagrangian fibration and $\mu$ is birational. Hence
$\pi$ is a rational Lagrangian fibration. By Proposition \ref{H2-for-Xdd},
$\Pic (\Xdd)$ is of rank 2 and the only primitive effective classes with
square zero are $h\pm me$, so $\pi^*[\OOO_{\PP^d}(1)]=h\pm me$.
By construction, $\pi^*[\OOO_{\PP^d}(1)]$ is represented by a divisor of the form
$$
D_t=\{\xi\in\Xdd\mid h^0(\III_\xi\otimes\EEE_t)\neq 0\}\subset\Xdd
$$
for generic $t\in M$, where $\EEE_t$ denotes a stable vector bundle on $X$
representing $t$.
The class $h+me$ has negative intersection with the generic fiber
$\PP^1$ of the Hilbert--Chow map $\Xdd\lra X^{(d)}$, hence has the whole
exceptional divisor $E$ in its base locus.
Hence, to see that $D_t\sim h - me$, it suffices to verify that $D_t$
does not contain $E$ as a fixed component. 
The support $\Supp \xi$ of a generic $\xi\in E$ is a set
of $d-1$ generic points of $X$. By Lemma \ref{technical}, (iii), the scheme
of zeros of any nonzero section $\sigma$ of $\EEE$ vanishing on $\Supp \xi$
is exactly $\Supp \xi$. But $\Supp\xi\subsetneqq\xi $, 
so $h^0(\III_\xi\otimes\EEE_t)= 0$, and $D_t\sim h - me$.
By the same argument as in Lemma \ref{fminuse}, $\pi$ is given by the complete
linear system $|h-me|$.
 \end{proof}

\begin{remark} \label{specialK3}
Though, as we mentioned in the introduction, $\pi$ is regular for generic $X$, it
may be nonregular for some special K3 surfaces with $\Pic X\simeq\ZZ^2$
which have a divisor class of degree $(2d-2)m^2$.
The following Proposition provides such a special K3 surface.
The map $\phi_{|h-me|}$ for this K3 surface is not regular, but a small deformation of $X$
kills its indeterminacy.
\end{remark}

\begin{proposition}\label{lattice-polar}
Let $X$ be a generic lattice-polarized K3 surface with
Picard lattice
$$
Q=\left(\begin{array}{cc}
2d & 2d-1+m\\
2d-1+m & 2d-2
\end{array}\right) .
$$
Let $f_{2d}$, $f_{2d-2}$ be effective classes forming a basis
of the Picard lattice in which the intersection form
is given by the above matrix.
Then the following properties are verified:

(i) If $d\geq 3$, then the linear system $|f_{2d}|$, resp. $|f_{2d-2}|$ embeds
$X$ into $\PP^{d+1}$, resp. $\PP^{d}$. If $d=2$, then $|f_{2d}|$ embeds $X$
as a smooth quartic in $\PP^3$, and $|f_{2d-2}|$ defines a double covering of
$\PP^2$.

(ii) Every curve in the linear systems $|f_{2d-2}|$, $|f_{2d}|$
is reduced and irreducible.

(iii) In addition to the rational map $\theta$ introduced in
(\ref{theta}), define the birational involution
$$
\iota :\Xdd\lrdash\Xdd\ ,\ \ \ \xi \mapsto (\langle \xi\rangle_{\PP^{d+1}}\cap X)-\xi\ ,
$$
where $\langle \xi\rangle_{\PP^{d+1}}$ denotes the linear span of 
$\xi$ in its embedding into $\PP^{d+1}$ by the linear system $|f_{2d}|$.
The corresponding isometry of the Bogomolov--Beauville lattice
is the reflection with respect to the vector $f_{2d}-e$ with square 2:
\begin{equation}\label{refl}
\iota^*:H^2(\Xdd)\lra H^2(\Xdd)\ ,\ \ \ c\mapsto -c+(c,f_{2d}-e)(f_{2d}-e)\ .
\end{equation}

Then the composite map $\pi =\theta\circ\iota$ is given by the complete
linear system $|f_{(2d-2)m^2}-me|$, where $f_{(2d-2)m^2}=
(m+1)f_{2d}-f_{2d-2}$
is an effective divisor class of degree $(2d-2)m^2$.
\end{proposition}

\begin{proof} This is similar to the work of Hassett--Tschinkel \cite{HTsch-1}
who produce on a K3  surface $X$ with two polarizations of degrees 4 and 8
an infinite series of polarizations $f_{2m^2}$ of degree $2m^2$ ($m\geq 2$)
and the abelian fibration maps on $X^{[2]}$ given by the linear system
$f_{2m^2}-me$.
The assertions (i), (ii) follow easily from the 
surjectivity of the period mapping
for K3 surfaces \cite{LP} and from the results
of \cite{SD}, \cite{Kov}. For formula (\ref{refl}) of part (iii),
see \cite{O'G-2}, 4.1.2. 
A direct calculation using (\ref{refl}) shows that
$f_{(2d-2)m^2}-me=\iota^*(f_{2d-2}-e)$. 
The class $f_{(2d-2)m^2}$ has positive square
and positive scalar product with $f_{2d}$, hence is effective.
As $\iota$ is an isomorphism in codimension 1, the dimensions of the linear
systems $|f_{2d-2}-e|$ and $|\iota^*(f_{2d-2}-e)|$ are the same, so
$\pi$ is defined by a complete linear system.
\end{proof}



\begin{thebibliography}{-----------}

\bibitem[ACGH]{ACGH} E. Arbarello, M. Cornalba, P. A. Griffiths,
J. Harris, {\em  
Geometry of algebraic curves, Vol. I}, 
Grundlehren der Mathematischen Wissenschaften,
Springer-Verlag, New York, 1985.


\bibitem[Beau-1]{Beau-1}  Beauville, A.: {\em
Some remarks on K\"ahler manifolds with $c_1 = 0$}, 
Classification of algebraic and analytic manifolds (Katata, 1982), 
1--26, Progr. Math., 39, Birkh\"auser, Boston, MA, 1983.


\bibitem[Beau-2]{Beau-2} Beauville, A.: {\em 
Vari\'et\'es K\"ahleriennes dont la premi\`ere classe de Chern est nulle},
J. Differential Geom. {\bf 18},  755--782 (1983).

\bibitem[Beau-3]{Beau-3} Beauville, A.: {\em 
Syst\`emes hamiltoniens compl\`etement int\'egrables associ\'es aux surfaces $K3$}, 
Problems in the theory of surfaces and their classification (Cortona, 1988), 25--31, 
Sympos. Math., XXXII, Academic Press, London, 1991.

\bibitem[Bo]{Bo} Bogomolov, F.A.: {\em 
Hamiltonian K\"ahler manifolds}, 
Sov. Math., Dokl. {\bf 19}, 1462-1465 (1978).

\bibitem[Bou]{Bou} Boucksom, S.: {\em 
Le c\^one k\"ahl\'erien d'une vari\'et\'e hyperk\"ahl\'erienne},
C. R. Acad. Sci., Paris, Ser. I, Math. {\bf 333}, 935-938 (2001). 

\bibitem[Cal-1]{Cal-1} C\u ald\u araru, A. {\em 
Derived Categories of Twisted Sheaves on Calabi-Yau Manifolds}, 
Ph.D. Thesis, Cornell University, 2000.

\bibitem[Cal-2]{Cal-2} C\u ald\u araru, A. {\em 
Nonfine moduli spaces of sheaves on $K3$ surfaces},  
Int. Math. Res. Not. {\bf 2002},  1027--1056 (2002).


\bibitem[F]{F}  Fujiki, A.: {\em 
On primitively symplectic compact Kahler V-manifolds of dimension four}, In: 
Classification of algebraic and analytic manifolds, Proc. Symp., Katata/Jap. 1982, Prog. Math. {\bf 39}, 71-250 (1983). 

\bibitem[GH]{GH} Griffiths, P. A., Harris, J.: {\em Principles of
Algebraic Geometry}, John Wiley \&\ Sons, New York, 1978.

\bibitem[Gu]{Gu} Gulbrandsen, M. G.: {\em Lagrangian fibrations on generalized Kummer varieties},
math.AG/0510145. 


\bibitem[Har-1]{Har-1} Hartshorne, R.: {\em Residues and Duality},
Lecture Notes in Math., No. 20,
Springer, 1966.

\bibitem[Har-2]{Har-2} Hartshorne, R.: {\em Algebraic geometry},
Graduate Texts in Mathematics, No. 52, Springer, 1977.


\bibitem[HasTsch-1]{HTsch-1} Hassett, B., Tschinkel, Yu.: {\em Abelian 
fibrations and rational points on symmetric products},  
Internat. J. Math. {\bf 11}, 1163--1176 (2000).

\bibitem[HasTsch-2]{HTsch-2} Hassett, B., Tschinkel, Yu.:
{\em Rational curves on holomorphic symplectic fourfolds},  
Geom. Funct. Anal. {\bf 11}, 1201--1228 (2001). 

\bibitem[Hu-0]{Hu-0} Huybrechts, D.:
{\em  
Birational symplectic manifolds and their deformations},
J. Differ. Geom. {\bf 45}, 488-513 (1997).

\bibitem[Hu-1]{Hu-1} Huybrechts, D.:
{\em Compact hyper-K\"ahler manifolds: basic results},
Invent. Math. {\bf 135}, 63--113 (1999); {\em Erratum}, ibid,
{\bf 152}, 209--212 (2003). 

\bibitem[Hu-2]{Hu-2} Huybrechts, D.:
{\em  The K\"ahler cone of a compact hyperk\"ahler manifold,} 
Math. Ann.  {\bf 326}, 499--513  (2003).



\bibitem[HL]{HL} Huybrechts, D., Lehn, M.:
{\em The Geometry of Moduli Spaces of Sheaves}, Aspects of Math.,
Vol. E 31, Friedr. Vieweg \&\ Sohn,
Braunschweig (1997). 

\bibitem[HS]{HS} Huybrechts, D., Stellari, P.:
{\em Equivalences of twisted K3 surfaces,} math.AG/0409030.


\bibitem[IR]{IR}  Iliev, A.,  Ranestad, K.: {\em
The abelian fibration on the Hilbert cube of a K3 surface of genus 9},
e-print math.AG/0507016.




\bibitem[Kov]{Kov}  Kov\'acs, S. J.: {\em The cone of curves of a
$K3$ surface},  Math. Ann.  {\bf 300}, 681--691 (1994).



\bibitem[L]{L} Lazarsfeld, R.: {\em Brill-Noether-Petri without degenerations},
J. Differential Geom. {\ bf 23}, 299--307 (1986). 

\bibitem[LeP-1]{LeP-1} Le Potier, J.:
{\em Sur l'espace de modules des fibr\'es de Yang et Mills,}
Mathematics and physics (Paris, 1979/1982), 65--137, Progr. Math., 37,
Birkh\"auser, Boston, 1983.


\bibitem[LeP-2]{LeP-2} Le Potier, J.: {\em
Faisceaux semi-stables et syst\`emes coh\'erents},
Vector bundles in algebraic geometry
(Durham, 1993), Ed. N. Hitchin et al., 179--239,
London Math. Soc. Lecture Note Ser., 208, 
Cambridge Univ. Press, Cambridge, 1995.


\bibitem[LP]{LP} Looijenga, E., Peters, C.: {\em  Torelli theorems 
for K\"ahler $K3$ surfaces},  Compositio Math. 
{\bf 42}, 145--186  (1980/81).

\bibitem[Mac]{Mac} Maclane, S.: {\em
Homology,} Springer, Berlin, 1963.

\bibitem[Mar]{Mar} Markman, E.: {\em 
Brill-Noether duality for moduli spaces of sheaves on 
K3 surfaces}, J. Algebraic 
Geometry {\bf 10}, 623--694 (2002).

\bibitem[Mat-1]{Mat-1} Matsushita, D.: {\em 
On fibre space structures of a projective irreducible
symplectic manifold},  Topology {\bf  38}, 79--83  (1999);
Addendum, ibid, {\bf  40}, 431--432  (2001).

\bibitem[Mat-2]{Mat-2} Matsushita, D.: {\em Higher direct images of 
dualizing sheaves of Lagrangian fibrations},  
Amer. J. Math. {\bf  127},  243--259  (2005).


\bibitem[Mor]{Mor} Morrison, D.: {\em
The geometry of K3 surfaces,} Lectures delivered at the Scuola 
Matematica Interuniversitaria, Cortona, 1988.


\bibitem[Mu-1]{Mu-1}  Mukai, S.: {\em 
Symplectic structure of the moduli space of sheaves 
on an abelian or $K3$ surface},  Invent. Math. {\bf  77},  101--116  (1984).

\bibitem[Mu-2]{Mu-2}  Mukai, S.: {\em On the moduli space of 
bundles on $K3$ surfaces. I}, Vector bundles on algebraic varieties 
(Bombay, 1984), 341--413, Tata Inst.
Fund. Res. Stud. Math., 11, Bombay, 1987. 


\bibitem[Mu-3]{Mu-3} Mukai,S.:
{\em  Duality of polarized $K3$ surfaces}, In:  New trends in algebraic
 geometry (Warwick, 1996),  311--326, London Math. Soc. Lecture Note Ser.,
 264, Cambridge Univ. Press, Cambridge, 1999.


\bibitem[O'G-1]{O'G-1} O'Grady, K. G.:
{\em The weight-two Hodge structure of moduli spaces of sheaves
 on a $K3$ surface},  J. Algebraic Geom. {\bf  6}, 599--644  (1997).

\bibitem[O'G-2]{O'G-2} O'Grady, K. G.: {\em 
Involutions and linear systems on holomorphic symplectic manifolds},
math.AG/0403519.

\bibitem[S-1]{S-1} Sawon, J.: {\em 
Abelian fibred holomorphic symplectic manifolds},
Turk. J. Math. {\bf 27}, 197-230 (2003). 

\bibitem[S-2]{S-2} Sawon, J.: {\em 
Lagrangian fibrations on Hilbert schemes of points on K3 surfaces},
e-print math.AG/0509224.


\bibitem[SD]{SD} Saint-Donat, B.: {\em Projective models of $K$-$3$ surfaces}, 
Amer. J. Math. {\bf 96}, 602--639 (1974).



\bibitem[Sim]{Sim} Simpson, C. T.: {\em Moduli of representations of the
fundamental group of a smooth projective variety I,} Publ. Math.
I.H.E.S. {\bf 79}, 47--129 (1994).


\bibitem[Tyu-1]{Tyu-1}  Tyurin, A. N.: {\em Special $0$-cycles on a polarized 
surface of type $K3$} (Russian), Izv. Akad. Nauk SSSR Ser. Mat. {\bf 51},
131--151 (1987).

\bibitem[Tyu-2]{Tyu-2}  Tyurin, A. N.: {\em Cycles, curves and vector 
bundles on an algebraic surface}, Duke Math. J. {\bf 54}, 1--26 (1987).



\bibitem[Y-1]{Y-1} Yoshioka, K.: {\em Moduli spaces of stable sheaves on abelian surfaces},  Math. Ann. {\bf 321}, 817--884  (2001). 

\bibitem[Y-2]{Y-2} Yoshioka, K.: {\em Moduli spaces of twisted sheaves on a projective variety},
math.AG/0411538

\end{thebibliography}
\end{document}